\input amstex
\documentstyle{gen-p}
\overfullrule=0pt      \input epsf 
  
\define\V{{\Cal V}}\define\A{{\Cal A}}

\input amssym.def
\define\bbZ{{\Bbb Z}} \define\bbR{{\Bbb R}} \define\bbQ{{\Bbb Q}}
\define\bbC{{\Bbb C}} \define\bbP{{\Bbb P}}  \define\g{{\frak g}}
\define\Hp{{\Bbb H}}

\topmatter
\title Zhu's algebra, the $C_2$ algebra, and twisted modules \endtitle
\author{Matthias R. Gaberdiel and Terry Gannon}\endauthor
\address{Institut f\"ur Theoretische Physik, ETH Z\"urich, CH-8093 Z\"urich, Switzerland}\endaddress \address{Math Department, University of Alberta, 
Edmonton, Canada T6G 2G1}\endaddress
\rightheadtext{Zhu's algebra, $C_2$, and twisted modules}
\email gaberdiel\@itp.phys.ethz.ch  tgannon\@math.ualberta.ca\endemail
\abstract 
In his landmark paper, Zhu associated two associative algebras to a VOA:
what are now called Zhu's algebra and the $C_2$-algebra. The former has 
a nice interpretation in terms of the representation theory of the VOA, while the 
latter only serves as a finiteness condition. In this paper we undertake first steps
to unravel the interpretation of the $C_2$ vector space. In particular we 
suggest that it sees and controls the twisted representations of the VOA.
\endabstract
\thanks This paper is based on work done during an extended visit by TG at ETH
funded by the CTS  -- he thanks them warmly for their generous hospitality.
The research of MRG is supported in part by the Swiss National Science Foundation,
while that of TG is supported in part by NSERC. We thank Andy Neitzke, who
collaborated in an early stage of this work.\endthanks

\endtopmatter

\document

\head I. Introduction\endhead

We study the modules of a finite group $G$ through an associative
algebra $\bbC G$ (its group algebra) with identical modules,
and we study the modules of a Lie algebra ${\frak g}$ through an
associative algebra $U{\frak g}$ (its universal enveloping algebra), again
with the same modules.
Likewise, the modules of a vertex operator algebra (VOA)
${\Cal V}$ are in natural one-to-one correspondence with those of an
associative algebra $\A(\V)$ called {\it Zhu's algebra}.

The most privileged class of Lie algebras are the finite-dimensional semi-simple
ones. The analogous notion for VOAs is termed
{\it rationality}: any sufficiently nice module of a rational VOA is
completely reducible into a direct sum of simple modules. For a Lie algebra
${\frak g}$, semi-simplicity is equivalent to the structural property that
the radical of ${\frak g}$ vanishes. For a VOA $\V$, the structural 
characterisation of rationality is conjectured to be the
finite-dimensionality and semi-simplicity of the algebra $\A({\Cal V})$. 

Closely related to Zhu's algebra is the $C_2${\it 
-algebra} $\A_{[2]}(\V)=\V/C_2(\V)$. The hypothesis of its
finite-dimensionality was first introduced and heavily
used in [{\bf 20}]. It implies for instance that $\V$ and its modules are finitely
generated, that $\V$ has only finitely many irreducible modules,
that the fusion coefficients are finite, and that the $\V$-modules $M$ have 
characters $\chi_M(\tau)$ holomorphic in $\Hp$ and obeying a weak modular 
invariance.
It had been often conjectured that the finite-dimensionality of $A_{[2]}(\V)$
is equivalent to rationality of $\V$, but a counterexample is the
so-called triplet algebra [{\bf 12}], and now we understand $A_{[2]}(\V)$
finite-dimensionality to characterise the {\it `finite logarithmic'} VOAs [{\bf 18}].
To prove the finite-dimensionality of $\A_{[2]}(\V)$ from the
semi-simplicity and finite-dimensionality of $\A(\V)$ is a fundamental problem
of VOA theory (first conjectured in [{\bf 20}]), and one of the motivations for 
our work.

For the remainder of this paper, let
${\V}$ be a {rational VOA}. For concreteness, we can  take this to mean:
\roster
\item{} $\V$ is a {\it simple} VOA (so ${\Cal V}$ is an irreducible 
${\Cal V}$-module for itself);

\item{}  as a ${\Cal V}$-module, ${\Cal V}$ is isomorphic to 
its contragredient;

\item{} $\V$ is of {\it CFT-type}: ${\Cal V}_n=0$ for $n<0$, ${\Cal V}_0=
\bbC |0\rangle$; and

\item{} $\V$ is {\it regular}: every  weak ${\Cal V}$-module is completely reducible\endroster
(these hypotheses imply then semi-simplicity of $\A(\V)$ as well as 
finite-dimensionality of $\A_{[2]}(\V)$ [{\bf 17}]). Under these conditions,
Zhu's theorem [{\bf 20}] applies, and the $\bbC$-span of the characters 
corresponding to the irreducible modules is SL$_2(\bbZ)$-invariant.
Also, Verlinde's formula holds, and the ${\Cal V}$-modules form a modular
tensor category [{\bf 15}].

In this paper we shall mainly deal with very specific examples: most
importantly, ${\Cal V}$ associated to
even lattices, and to affine algebras at integral level.
For the construction and basic properties of these VOAs, see 
{\it e.g.}\ [{\bf 16}];
for other general surveys of VOAs {\it etc.}\ see 
{\it e.g.}\ [{\bf 5,13}].

The authors have been exploring the intimate relation between Zhu's algebra
and the $C_2$-algebra in rational VOAs [{\bf 11}]; in the present paper
we motivate and review some initial results and early speculations. 
We describe in Section II Zhu's algebra, the $C_2$-algebra,
and the notion of twisted modules. Zhu's definitions [{\bf 20}] of $\A(\V)$
and $\A_{[2]}(\V)$ are very technical and unmotivated; we include a motivation
due to Gaberdiel and Goddard  [{\bf 10}] which underlies our approach. We review
how Zhu's algebra sees the ${\V}$-modules, and explain how the closely related
$\A_{[2]}({\V})$ sees the twisted $\V$-modules.
Section III presents some of the general theory: {\it e.g.}\ 
the dual space $\A(\V)^*$  embeds naturally
in $\A_{[2]}({\Cal V})^*$; what meaning can be ascribed to the `discrepancy'
$\A_{[2]}({\Cal V})^*/\A(\V)^*$ and how often is it zero?  Section IV summarises
some of our calculations. We collect in Section V some speculations and
open questions.

It is an honour to dedicate this paper to Geoff Mason -- his work has been
an inspiration to both of us.

\head II. Background\endhead

\subhead II.1. Zhu's algebra: the what and the why\endsubhead
Let us begin by describing informally the ideas behind Zhu's algebra (surely 
at least some of this fed Zhu's intuition), because
it is crucial to the rest of this paper. As with much of VOA theory, the
basic ideas come from conformal field theory (CFT) or, equivalently,
perturbative string theory.

The correlation functions are how a quantum field theory
makes contact with experiment. All physical content is there. In the case
of a CFT, the correlation functions are (essentially) the conformal blocks.
Consider the Riemann sphere $\bbP^1$, and choose $n\ge 2$ distinct points
$w_i\in\bbP^1$; to each of these points we select an irreducible ${\Cal V}$-module
$M^i$, and a vector $a_i\in M^i$. The corresponding conformal blocks will 
be a span of formal objects of the form
$$\langle {\Cal Y}_1\cdots {\Cal Y}_{n-2}\rangle\ ,\eqno(2.1)$$
where each ${\Cal Y}_i$ is an intertwining operator -- {\it e.g.}\ ${\Cal Y}_1$
will be a linear map sending $a_1\otimes a_2\in M^1\otimes M^2$ to a formal
series with coefficients in some irreducible ${\Cal V}$-module $N^1$ contained
in the fusion product of $M^1$ and $M^2$. 
The details are not important:
the conformal block will be a complex-valued function, multilinear in
the $a_i\in M^i$, and locally meromorphic (as a density or differential form) in 
the configuration space $(w_1,\ldots,w_n)\in (\bbP^1)^n\setminus 
\Delta$, where $\Delta$ consists of the diagonals, where some $w_i$ coincides with
some $w_j$.
The space of all such conformal
blocks, has dimension given by a Verlinde formula and depends only on $M^i$.
It suffices to consider only {\it highest-weight states} $a_i\in (M^i)_0$,
{\it i.e.}\ the vectors in $M^i$ of lowest $L_0$-eigenvalue, as the other conformal
blocks are obtained from these by standard differential operators.

Let us probe the (vacuum-to-vacuum) conformal block (2.1) by `inserting' a state
$v\in{\Cal V}$.
This amounts to adding an ($n+1$)-st point $w_{n+1}$, with module $M^{n+1}={\Cal V}$.
Conventionally we shall take, without loss of generality, $w_{n+1}=0$ (which we think 
of as time $t=-\infty$), so the  other $w_i$ are now required to avoid 0. 
The resulting conformal block is a function
$$\eta_{\{(w_1,M^1),\ldots,(w_n,M^n)\}}\in\text{Hom}\Bigl({\Cal V},
(M^1\otimes\cdots \otimes M^n)^*\Bigr)\  , \eqno(2.2)$$
where `$\otimes$' denotes fusion product, `$*$' denotes a restricted dual,
and `Hom' is the space of ${\Cal V}$-module maps (intertwiners).
This looks fancier than it really is: `Hom' chooses the specific conformal
block (from among the space of these); the conformal block is a complex-valued
function of $v\in{\Cal V}$ and the $a_i\in M^i$. Choosing $v$ to be the vacuum
$|0\rangle$ in the VOA recovers the conformal block (2.1).

There is a large subspace ${\Cal O}_{(w_1,\ldots,w_n)}$ of ${\Cal V}$,
independent of the choice of modules $M^i$, 
for which (2.2) must vanish for elementary reasons. For example, if say
$w_1=\infty$ (and the $a_i$ are indeed highest-weight), then 
${\Cal O}_{(w_1,\ldots,w_n)}$ is spanned by 
$$\text{Res}_z\left(Y(a,z)z^{-1-m+(2-n)|a|}\prod_{j=2}^n(z-w_j)^{|a|}b\right)\eqno(2.3)$$
for all choices of integer $m>0$, and vectors $a,b\in{\Cal V}$, where $a $ is
homogeneous with respect to $L_0$: $L_0a=|a|\,a$.
Indeed, such a vector can be (formally) written as  a
contour integral, and it is easy to see that the integrand does not have
any poles apart from at 0.

Let us define $\A_{(w_1,\ldots,w_n)}({\Cal V})={\Cal V}/{\Cal O}_{(w_1,\ldots,w_n)}$.
Then for any  $a_i\in (M^i)_0$,  we can regard
$\eta_{\{(w_1,M^1),\ldots,(w_n,M^n)\}}$ as a linear functional
in $\A_{(w_1,\ldots,w_n)}^*$. The converse is also true but much deeper:
any $\eta\in \A_{(w_1,\ldots,w_n)}^*$ is a correlation function, for some
choice of modules $M^i$ and highest weight states $a_i\in (M^i)_0$ attached to 
the points $w_i\in \bbP^1$. For a rational VOA, all these spaces will be 
finite-dimensional.

An elementary observation [{\bf 19}] is that this space 
$\A_{(w_1,\ldots,w_n)}({\Cal V})$
is independent of the order of the $w_i$'s. Moreover,
by an analytic continuation argument based on a Knizhnik--Zamolodchikov-like 
connection, for each homotopy class of paths in 
the configuration space $(\bbP^1\setminus 0)^{n}\setminus\Delta$ linking two 
$n$-tuples $(w_1,\ldots,w_n)$ and $(w'_1,\ldots,w_n')$, we get a natural
isomorphism $\A_{(w_1,\ldots,w_n)}({\Cal V})\cong \A_{(w_1',\ldots,w_n')}({\Cal V})$.
This means we get an action of the braid group $B_n$ on $\A_{(w_1,\ldots,w_n)}({\Cal V})$ by automorphisms, although this is not important for us.

Comparing (2.3) with [{\bf 20}], we see that
Zhu's algebra $\A(\V)$ is $\A_{(\infty,-1)}({\Cal V})$. In this case, the correlation functions
will vanish unless $M^1$ and $M^2$ are duals (contragredients) of each other.
So what we find from the above is that Zhu's algebra carries information
on all the irreducible ${\Cal V}$-modules: its dual
$\A({\Cal V})^*$ naturally decomposes into a sum of nontrivial subspaces,
one for each pair $(M,M^*)$. We will make this more precise shortly.

The algebraic structure of Zhu's algebra comes from the action of the 
{\it zero-modes} $o(a)=a_{|a|-1}$ on the highest weight states. Indeed, 
$o(a)$ commutes with $L_0$ and thus maps each homogeneous 
subspace $M_n$ of any (irreducible) $\V$-module
$M=\oplus_{n\in\bbZ_{\ge 0}}M_n$ to itself, in particular the highest
weight space $M_0$.  We can describe
this zero mode action in terms of a product in the $(n+1)$-st  module 
$M^{n+1}={\Cal V}$; for example, for the case considered by Zhu 
where we just have two points $w_1=\infty$ and $w_2=-1$, the relevant
product is [{\bf 20}] 
$$a*b= \text{Res}_z\left(Y(a,z)z^{-1} (z+1)^{|a|}b
\right)\eqno(2.4)$$
for any $a,b\in\V$. By similar contour integral arguments as above it is easy to see
$$
\eta_{\{(\infty,M^1),(-1,M^2)\}}(a*b) [m_1 \otimes m_2] = 
\eta_{\{(\infty,M^1),(-1,M^2)\}}(b) [o(a) m_1 \otimes m_2]  \ , \eqno(2.5)
$$
where $m_i\in (M^i)_0$, $i=1,2$. As is expected from 
(2.5)  ${\Cal O}({\Cal V}) = {\Cal O}_{(\infty,-1)}({\Cal V})$ is an ideal 
for this action. Thus we obtain an action of ${\Cal V}$ on ${\Cal A}({\Cal V})$, and 
this turns ${\Cal A}({\Cal V})$  into an associative algebra. In fact, the product
structure simply describes the product of the zero modes: {\it i.e.}\ 
on highest weight states one finds $o(a*b) = o(a) o(b)$.

In general, the other spaces $\A_{(w_1,\ldots,w_n)}(\V)$ do not seem to 
naturally be algebras. However, their duals
 $\A_{(w_1,\ldots,w_n)}(\V)^*$ always have $n$ commuting actions
of $\A(\V)$ (one for each point $w_i$) -- see [{\bf 9}] for details.

For rational VOAs, Zhu [{\bf 20}] proved that his associative algebra $\A(\V)$ is 
finite-dimensional and semi-simple, and hence by Wedderburn's Theorem a direct sum 
of matrix algebras: in fact we have
$$\A(\V)\cong \oplus_M \text{End}(M_0)=\oplus_M M_0^*\otimes M_0\ ,\eqno(2.6)$$
where $M$ runs over all irreducible $\V$-modules.
 This is the precise sense in which Zhu's algebra
sees all $\V$-modules. The  2-point correlation
functions are para\-me\-trised by a choice of module $M$ and states $a_1\in M_0$ and
$a_2\in(M^*)_0=(M_0)^*$, together filling out the summand $M_0\otimes M_0^*$
in the dual of (2.6).

\subhead II.2. When 2 become 1: the $C_2$ algebra \endsubhead
The conformal blocks (2.1) live on the moduli space of an $n$-punctured sphere.
For example, for $n=4$, this moduli space can be identified with the sphere
minus 3 points (the M\"obius transformations send the 3 points to 0, 1, $\infty$, 
and the fourth point is then free provided it avoids those). Now in
string theory, it is meaningful (`amplitude factorisation') to move towards
boundary points
in these moduli spaces, {\it e.g.}\ to send two of these $n$ points together.
This is also meaningful in number theory (`cusps'), and the Deligne-Mumford
compactification of moduli space tells us to interpret those boundary points
as surfaces with nodes. For example, for $n=4$, the three boundary components
correspond to the three ways to partition the 4 points into 2 pairs; the
corresponding surface is two tangential spheres, each containing two of the 
points.

In any case, we are led to consider what happens when some of the $w_i$
coincide. Most of the treatment of the previous subsection goes through
without change. We can formally speak of conformal blocks, and define
the subspaces ${\Cal O}_{(w_1,\ldots,w_n)}(\V)$ and the corresponding
quotient spaces $A_{(w_1,\ldots, w_n)}(\V)$ as before. To understand
more explicitly what happens consider Zhu's algebra $\A(\V)$ which is 
naturally  isomorphic to  $\A_{(\infty,w_2)}(\V)$ for any $w_2\ne\infty$. For any
such $w_2$, the elements spanning ${\Cal O}_{(\infty,w_2)}(\V)$  are of the form
(after appropriate rescaling)
$$ \text{Res}_z \left(Y(a,z)z^{-2}\bigl(1-\frac{z}{w_2}\bigr)^{|a|}b\right)
= \sum_{l=0}^{|a|} {|a| \choose l} \, \bigl(-w_2\bigr)^{-l} \, 
a_{-2+l} \, b\ .
$$
In the limit $w_2\rightarrow \infty$, only the leading term survives, and thus
the space ${\Cal O}_{(\infty,\infty)}(\V)$ is spanned by the elements of the form
$a_{-2} \, b$; this is precisely  the $C_2({\Cal V})$ space of Zhu, and hence
$\A_{(\infty,\infty)}(\V)$ is the $C_2$ quotient space of Zhu. Similarly, 
the product  $a*b$ in $\A_{(\infty,w_2)}(\V)$ (appropriately rescaled) can be 
written as
$$a*_{w_2}b=\text{Res}_z \left(Y(a,z)z^{-1}\bigl(1-\frac{z}{w_2}\bigr)^{|a|}b\right)
= \sum_{l=0}^{|a|} {|a| \choose l} \, \bigl(-w_2\bigr)^{-l} \, 
a_{-1+l} \, b\ .
$$
This tends to the product $ab:= a_{-1}\, b$ of Zhu for the $C_2$-space 
$\A_{(\infty,\infty)}(\V)$. The resulting algebra is both commutative and associative, 
for elementary  reasons. In fact, Zhu noticed that this quotient  also has
a Poisson structure: the Lie bracket can be defined by $\{a,b\}:=a_0b$.
We shall call this (commutative Poisson) algebra the $C_2$-algebra 
$\A_{[2]}(\V)$. 


We can thus think of Zhu's algebra as being a deformation of the 
$C_2$-algebra. However, this deformation picture is certainly false 
if taken too literally (as it has in the literature
-- see {\it e.g.}\ the last paragraph of Section 4 of [{\bf 1}]).
In particular, although the ideal ${\Cal O}_{(\infty,w_2)}(\V)$ tends to
${\Cal O}_{(\infty,\infty)}(\V)$, elements that were non-trivial in the 
former can tend to 0 in the limit. The correct statement is that there is
a natural surjection 
${\Cal O}_{(\infty,w_2)}(\V)\rightarrow {\Cal O}_{(\infty,\infty)}(\V)$ (which 
may not be an injection). Thus the dim$\,\A(\V)$ of Zhu's algebra is bounded 
above by that of the $C_2$-algebra $\A_{[2]}(\V)$, and the dual $\A(\V)^*$ can
be regarded as a subspace of $\A_{[2]}(\V)^*$. We return to this in Section
III.1.

As before we can reorder the $n$ points $w_i$ so that the first $n_1$ are
identical, the next $n_2$ are identical (but different from the first $n_1$),
{\it etc}, where $n_1\ge n_2\ge n_k>0$ is some partition of $n$; then the
analytic continuation argument shows that the spaces $\A_{(w_1,\ldots,w_n)}(\V)$
and $\A_{(w_1',\ldots,w_n')}(\V)$ are isomorphic whenever the $w_i$ and the
$w_j'$ correspond to the same partition of $n$. Thus we may speak of
the space $\A_{[n_1,\ldots,n_k]}(\V)$. This explains our notation $\A_{[2]}(\V)$
for the $C_2$-algebra; likewise, Zhu's algebra $\A(\V)$ is $\A_{[1,1]}(\V)$
in this notation.

It is elementary that Zhu's algebra $\A(\V)$ sees two commuting actions of
the automorphism group Aut$(\V)$ of the VOA, one attached to each point $w_i$.
As these points are brought together to form $\A_{[2]}(\V)$, what survives
is the diagonal action. So the $C_2$-algebra carries an adjoint action of
Aut$(\V)$, helping significantly to organise $\A_{[2]}(\V)$, which in
specific calculations can get quite large. For lattice VOAs $\V_L$, Aut$(\V_L)$
contains the automorphism group of the lattice $L$; for affine algebra VOAs
$\V_{\g,k}$, Aut$(\V_{\g,k})$ contains the simply connected
Lie group corresponding to $\g$.

The importance of Zhu's algebra is that its representation theory is
isomorphic to that of the VOA.
On the other hand, it is hard to imagine any useful direct
relation between the $\A_{[2]}(\V)$-modules and the $\A(\V)$- or $\V$-modules.
As an algebra, $\A_{[2]}(\V)$ is isomorphic to the $d\times d$
diagonal matrices, where $d=\text{dim}\,\A_{[2]}(\V)$. Hence there are
exactly $d$ irreducible $\A_{[2]}(\V)$-modules, all one-dimensional: the
$i$th one is the projection to the $i$th diagonal entry of the matrices.
Nevertheless, we will explain next subsection that $\A_{[2]}$, or rather its
dual space, {\it is} intimately connected to the representation theory of the
VOA.

\subhead II.3. Twisted modules for lazy people\endsubhead
We shall assume the reader is familiar with the usual notion of a VOA 
module -- see {\it e.g.}\  [{\bf 16,5,13}] for more details.
Twisted modules are a natural generalisation, and a central part of the whole
VOA story. Indeed, they are key to the orbifold construction.
They are at least as important for VOAs, as projective representations are
 to groups. In fact they are sort of a dual concept to projective representation:
to unprojectify a {\it projective} representation, you take a central {\it extension}
of the group; to untwist a {\it twisted} module, you restrict to a {\it subalgebra}
of the VOA.

Probably the easiest path to twisted modules is through the loop algebra.
Let ${\frak g}$ be a finite-dimensional simple Lie algebra (over ${\Bbb C}$).
By the {\it loop algebra} ${\Cal L}\g$ we mean the space of all combinations
$\sum_{n\in\bbZ} a_nt^n$, where $a_n\in\g$ and all but finitely many $a_n$ are
0 ($t$ is a formal variable). This inherits a Lie algebra structure from
$\g$. The nontwisted affine Kac-Moody algebra $\g^{(1)}$ is just the extension
of $\g$ by a central element $c$ and a derivation $\ell_0$.

Now let $\alpha$ be any automorphism of $\g$, of order $N<\infty$. We can
diagonalise $\alpha$: for $0\le j<N$ let $\g_j$ be the eigenspace of $\alpha$
in $\g$ with eigenvalue $\xi_N^j$, where we write $\xi_N=e^{2\pi\imath/N}$.
Of course $\alpha$ extends to
an automorphism of ${\Cal L}\g$ by sending $t^n$ to $\xi_N^nt^n$,
and to the affine algebra $\g^{(1)}$ by fixing $c$ and $\ell_0$. By the
{\it twisted affine algebra} $\g^{(N)}$ we mean the subalgebra of
$\g^{(1)}$ fixed by $\alpha$. The twisted affine algebras behave very
similarly to the more familiar nontwisted ones.

Now let $\rho$ be any integrable highest weight representation of $\g^{(N)}$.
We can lift $\rho$ to $\g^{(1)}$ by defining $\rho(at^n)=\xi_N^{j+n}
\rho(at^{-j})t^{j+n}$ for $a\in\g_j$. This will not be a true representation of
the nontwisted algebra $\g^{(1)}$, as it obeys
$$[\rho(at^n),\rho(bt^m)]=\xi_N^{j+k+m+n}\rho([at^n,bt^m])\eqno(2.7)
$$
when $a\in\g_j$ and $b\in\g_k$.
We call such a $\rho$ a {\it twisted representation} of 
$\g^{(1)}$.
Thus a true representation of a twisted affine algebra lifts to a twisted
representation of a nontwisted affine algebra.

The definition for VOAs is very analogous (see {\it e.g.}\ [{\bf 5,13}]).
Incidentally, it is possible to generalise the spaces 
$\A_{(w_1,\ldots,w_k]}({\Cal V})$ of Section II.1
to the case where now at some (or all) of the $w_i$ 
states from a twisted $\V$-module $M^i$ are inserted-- see {\it e.g.}\ 
[{\bf 3}].

Twisted modules are a crucial, though unexplored, part of the $C_2$-algebra
story. We explained
at the end of Section II.1 how, for any $\V$-module $M$, any choice $u\in M_0
,v\in M_0^*$
yields a unique vector $u\otimes v\in \A(\V)^*$. Since $\A(\V)^*$ embeds in
$\A_{[2]}(\V)^*$, $u\otimes v$ can also be regarded as a vector in $\A_{[2]}
(\V)^*$. If instead $M$ is a twisted $\V$-module, then $u\otimes v$ maps
into the appropriate twisted Zhu's algebra $\A_g(\V)$, defined in [{\bf 3}].
Implicit in the above treatment is that twisted modules are characterised
by monodromy properties about the point $w$ they have  been inserted; as the
two points $w_i$ are brought together, we cannot tell any more whether
$u\otimes v$ came from twisted or untwisted modules. This means that
each $\A_g(\V)^*$ also embeds into $\A_{[2]}(\V)^*$. Clearly, the images for
different automorphisms $g$ can overlap, and we do not yet understand the
relation between these different images. But it should be clear that
the $C_2$-algebra must be large enough to contain every $\A_g(\V)$.
This accounts for some, and perhaps all, instances where the $C_2$-algebra
is larger than Zhu's algebra. It also provides a partial, and perhaps complete,
answer to the question of the direct relevance of the $C_2$-algebra (or rather
its dual) to the representation theory of $\V$.

\head III. Abstract nonsense\endhead

In this section we collect some general comments about the $C_2$-algebra and
its relation with Zhu's algebra. 

\subhead III.1. Zhu's algebra as a deformation of the $C_2$-algebra\endsubhead
As we have explained before in Section II.2, Zhu's algebra $\A(\V)$ is a 
`deformation' of the $C_2$ algebra $\A_{[2]}(\V)$. 
As we have also explained there, the dimension of $\A(\V)$ may be smaller 
than that of the $C_2$ space $\A_{[2]}(\V)$. The situation is vaguely reminiscent of 
deformation quantisation, 
where a commutative Poisson algebra (describing the classical world) is deformed
into a noncommutative algebra (describing the quantum world). For this reason 
we suggest calling a VOA {\it anomalous} if the dimension of $\A_{[2]}(\V)$
is strictly larger than that of $\A(\V)$.

Note that  $\A(\V_1\otimes \V_2)=\A(\V_1)\otimes\A(\V_2)$ and 
$\A_{[2]}(\V_1\otimes \V_2)=\A_{[2]}(\V_1)\otimes\A_{[2]}(\V_2)$, so 
the $C_2$-algebra and Zhu's algebra of the tensor product $\V_1\otimes\V_2$
of VOAs will have equal dimension iff the same holds for both $\V_1$ and
$\V_2$.

As explained in Section II.2, we can think of the dual $\A(\V)^*$ as being a
subspace of $\A_{[2]}(\V)^*$. Let us call the quotient $\A_{[2]}(\V)^*/\A(\V)^*$
the {\it deficiency}, for want of a better name. This finite-dimensional space
is then nontrivial iff $\V$ is anomalous. Is there a cohomological interpretation
for the deficiency? Of course there is a rich relation of Hochschild cohomology
to the deformation theory of algebras [{\bf 14}]. For example, the
group $H^i(\A_{[2]}(\V);\A_{[2]}(\V))$ for $i=1,2$ respectively, equals the space
of infinitesimal automorphisms, respectively the space of infinitesimal deformations,
of the $C_2$-algebra, and this group for $i=3$ controls whether these 
infinitesimal deformations can be `integrated'. Hence 
whenever that second cohomology group vanishes, 
the VOA will either be anomalous, or dim$\,M_0=1$ for all irreducible $M$.
However, this remark is too naive to be of any value, because the
$C_2$-algebra is too uninteresting. A proper cohomological treatment of
deficiency {\it etc.}\ would have to involve more of the structure of $\V$.

\subhead III.2. Filtrations versus gradings\endsubhead
An algebra $A$ is called {\it graded} if $A$ is the direct sum $\oplus_{n=0}^\infty
A^n$ of subspaces $A^n$, such that $A^mA^n\subseteq A^{m+n}$. For example,
the polynomials $A=\bbC[x]$ are graded by degree, so each $A^n=\bbC x^n$ here
is one-dimensional.

An algebra $A$ is called {\it filtered} if $A$ is the union $\cup_{n=0}^\infty
A_n$ of an increasing sequence $A_0\subseteq A_1\subseteq A_2\subseteq\cdots$
of subspaces, such that $A_mA_n\subseteq A_{m+n}$. Any graded algebra is
filtered: just take $A_n=\oplus_{m=0}^n A^m$. A filtered algebra which is
not graded, is the universal enveloping algebra $U\g$: assign degree 1 to
every element of $\g$, and let $U\g_n$ consist of all polynomials in
$\g$, each term in which has total degree $\le n$. Degree does not define
a grading on $U\g$ (unless $\g$ is abelian): for any noncommuting $x,y\in\g$,
$xy$ and $yx$ both have degree 2 but their difference $[x,y]$ has degree 1.

The $C_2$-algebra is graded by $L_0$-eigenvalue, since its ideal $C_2(\V)$
is spanned by homogeneous elements $a_{-2}b$, and the product $a_{-1}b$
respects $L_0$-grading. On the other hand, Zhu's algebra is only filtered
by $L_0$, since the elements (2.3) spanning its ideal are not homogeneous.

There is a standard way to go from a filtered algebra $A=\cup_n A_n$ to a
graded algebra $A_{\text{gr}}$: define $(A_{\text{gr}})^n=A_n/A_{n-1}$. If $A$ is in fact
graded, then $A_{\text{gr}}\cong A$. If $A$ is finite-dimensional, then
dim$\,A_{\text{gr}}=\text{dim}\,A$. For example, $U\g_{\text{gr}}$ is
naturally isomorphic to the symmetric (polynomial) algebra $S\g$, obtained
by identifying $\g$ with $U\g_1/U\g_0$. $U\g$ carries two commuting
$\g$-actions: the left- and right-regular actions $gu$ and $-ug$; $S\g$
carries the adjoint $\g$-action $gu-ug$.

It is elementary to verify that the `gradification' $\A(\V)_{\text{gr}}$
can be identified (though not canonically) with a subspace of $\A_{[2]}(\V)$,
and hence with all of $\A_{[2]}(\V)$ if their dimensions match.
What role in the general story this gradification plays, is not yet clear to us.
But as we shall discuss in Section III.4, for the VOAs associated to
affine algebras, this point of view could be very important.

\subhead III.3. Zhu's algebra and the $C_2$-algebra for lattices\endsubhead
Let $L$ be any even positive-definite lattice (so $\alpha\cdot\alpha\in 2\bbZ_{\ge 0}$
for any $\alpha\in L$). Let $n$ be its dimension. Fix a basis $\{\beta_1,
\ldots,\beta_n\}$ of $L$. See {\it e.g.}\  [{\bf 16]} for the construction of $\V_L$.
As a vector space, $\V_L$ is spanned by terms of the form
$$\beta_{i_1}(-k_1)\,\beta_{i_2}(-k_2)\cdots\beta_{i_m}(-k_m)\,e^\alpha\ , \eqno(3.1)$$
where $m\ge 0$, each $k_i\in \bbZ_{>0}$, and $\alpha\in L$. The oscillators
$\beta_i(-k)$ commute with each other -- apart from that, the vectors in
(3.1) are linearly independent.

It can be shown [{\bf 11}] that the $C_2$-algebra ideal $C_2(\V_L)$ is spanned
by all terms of the form (3.1), provided at least one $k_i$ is $\ge 2$, together
with all vectors of the form
$$\beta_{i_1}(-1)\cdots\beta_{i_m}(-1)\,\gamma(-1)^{\text{max}\{0,1+\gamma\cdot
\gamma-|\gamma\cdot\alpha|\}}e^\alpha\, .\eqno(3.2)$$
Thus a basis for $\A_{[2]}(\V_L)$ can be found with coset representatives
of the form
$$\beta_{i_1}(-1)\cdots\beta_{i_m}(-1)\,e^\alpha \ ,\eqno(3.3)$$
where $\alpha$ belongs to the finite set 
$$S_L=\{\alpha\in L\,|\,\gamma\cdot\gamma\ge\gamma\cdot\alpha\ \forall\gamma\in
 L\}$$
of `small' lattice vectors; of course which oscillators $\beta_{i_j}(-1)$ 
to choose in (3.3) depends very much on the choice of $\alpha\in S_L$. 

It is easy to see from
this description of $C_2(\V_L)$ that $\A_{[2]}(\V_L)$
is finite-di\-men\-sio\-nal for any $L$ (first proved in [{\bf 4}]).
Next section we explain how to
use the preceding paragraph
 to find $\A_{[2]}(\V_L)$, or at least lower bounds for dim$\,\A_{[2]}(
\V_L)$, for explicit $L$.

The irreducible modules for $\V_L$ are in natural one-to-one correspondence 
with the cosets $[t]\in L^*/L$, where $L^*$ is the dual lattice of $L$. 
The character of the module corresponding to ${[t]}$ is the theta series of
the shifted lattice $[t]$, divided by $\eta(\tau)^n$. Its leading term is
the number $N_{[t]}$ of vectors in $[t]$ of smallest norm. The dimension
of Zhu's algebra is then 
$$\text{dim}\,\A(\V_L)=\sum_{[t]\in L^*/L}N_{[t]}^2\ .\eqno(3.4)$$
A priori, there seems little relation between (3.4) and dim$\,\A_{[2]}(\V_L)$
-- a reason for this is implicit in Section IV.4.

\subhead III.4. Affine Lie algebras\endsubhead
An important and nontrivial class of rational VOAs are associated to a
choice of finite-dimensional simple Lie algebra $\g$, and a positive integer
$k$ (the `level'). The associated rational VOA was constructed in [{\bf 8}]
and will be denoted $\V_{\g,k}$. Its homogeneous space $(\V_{\g,k})_1$
is canonically identified with $\g$. This VOA is intimately connected to the
affine nontwisted algebra $\g^{(1)}$; in particular, as spaces $\V_{\g,k}$
is the integrable $\g^{(1)}$-module $L(k\Lambda_0)$, and the irreducible
$\V_{\g,k}$-modules are the  level $k$ integrable highest weight 
$\g^{(1)}$-modules $L(\lambda)$. Write 
$\lambda=\sum_{i=0}^r\lambda_i\Lambda_i$.

Zhu's algebra here can be identified [{\bf 8}] with the quotient
$U\g/\langle e_\theta^{k+1}\rangle$, where $\langle e_\theta^{k+1}\rangle$
is the 2-sided ideal of $U\g$ generated by $e_\theta^{k+1}$ ($\theta$ is
the highest root of $\g$). The space $M_0$ for the $\V_{\g,k}$-module
associated to $\lambda$, can be identified with the irreducible $\g$-module
with highest weight $\overline{\lambda}=\sum_{i=1}^r\lambda_i\Lambda_i$,
so the dimension of Zhu's algebra then follows from {\it e.g.}\ 
Weyl's dimension formula.

The $C_2$-algebra arises naturally as a quotient $S\g/I(k)$. Here, $S\g$
is generated by the $-1$-modes of $(\V_{\g,k})_1\cong\g$, and the $\g$-action
on it comes from the zero-modes of $(\V_{\g,k})_1\cong\g$. The $m$'th graded
piece of $S\g$ can be identified with the $m$'th symmetric power of the
adjoint module of $\g$. The ideal $I(k)$ is generated from
$e_\theta^{k+1}$ using the  $\g$-action on $S\g$ described earlier.

Zhu's algebra inherits the filtration of $U\g$.
Put $I_n=\langle e_\theta^{k+1}\rangle\cap U\g_n$ and write $I_{\text{gr}}
=\oplus_{n} I_n/I_{n-1}$ as usual. Then the `gradation'
$\A(\V)_{\text{gr}}$ is canonically isomorphic to $S\g/I_{\text{gr}}$.
We would like to understand better the relation between the ideals
$I_{\text{gr}}$ and $I(k)$  of $S\g$, as this seems a very promising
approach to the question of anomalous $\V_{\g,k}$. The former ideal
contains the latter, and this defines the surjection $\A_{[2]}(\V_{\g,k})
\rightarrow \A(\V_{\g,k})$.
For most pairs $\g,k$ it seems, these ideals are identical (see Section IV.5
below).

\head IV. Calculations\endhead

\subhead IV.1. The Virasoro minimal models\endsubhead
Perhaps the easiest examples to work out are the Virasoro minimal models
$\V^{\text{Vir}}_{p,q}$, where $p,q\in\{2,3,4,\ldots\}$ are coprime (see
{\it e.g.}\ [{\bf 6}]).
In this case there are $(p-1)(q-1)/2$ irreducible modules $M$, all with
1-dimensional $M_0$, so $\A$ here is commutative, of dimension $(p-1)(q-1)/2$.
$\A_{[2]}$ is easy to identify because one null vector is $L_{-1}|0\rangle$,
so $\A_{[2]}$ has a basis of the form $L_{-2}^i|0\rangle+C_2$; the other null 
vector, whose leading term is $L_{-2}^{(p-1)(q-1)/2}|0\rangle$, then forces 
$0\le i<(p-1)(q-1)/2$. Thus the minimal models are non-anomalous.

\subhead IV.2. Affine $sl(2)$ at level $k$ \endsubhead
This is again very easy, and we know of at least 4 independent ways to
prove that the VOA is non-anomalous. For reasons of space we shall give only one.

Let $k$, the level, be any positive integer. The rational VOA $\V_{sl(2),k}$,
as a space, is given by the highest weight $sl(2)^{(1)}$-module $L_{k\Lambda_0}
=U(sl(2)^{(1)})|0\rangle$, and so inherits the filtration from the universal 
enveloping algebra. This permits us to refine the character
of $\V_{sl(2),k}$, to be a function not only of the usual $q$ (which keeps
track of the $L_0$-eigenvalue, what we are calling the grade) and $z$ (which
lies in the $SL(2)$ maximal torus so is the argument for $SL(2)$-characters),
but another parameter $t$ (which will keep track of this degree). More precisely,
each creation operator $x_{-n}$ will contribute 1 to the degree but $n$ to
the grade.

\noindent The result is [{\bf 7}]:
$$\chi_{\V_{sl(2),k}}(q,z,t)=\sum_{\vec{h},\vec{e},\vec{f}\in\bbZ_{\ge 0}^k}
t^{|\vec{e}|+|\vec{h}|+\vec{f}|}z^{2(|\vec{e}|-|\vec{f}|)}
\frac{q^{\vec{e}A\vec{e}^T+\vec{h}A\vec{h}^T+\vec{f}A\vec{f}^T
+\vec{e}B\vec{h}^T+\vec{h}B\vec{f}^T}}{(q)_{\vec{e}}(q)_{\vec{f}}(q)_{\vec{h}}} \ ,
\eqno(4.1)$$
where for $\vec{n}\in\bbZ_{\ge 0}^k$ we set $|\vec{n}|=\sum_{i=1}^k in_i$,
and $(q)_{\vec{n}}=\prod_i (q)_{n_i}$ where $(q)_n=\prod_{j=1}^n(1-q^j)$.
The $k\times k$ matrices $A$ and $B$ are defined by $A_{ij}=\text{min}\{i,j\}$
and $B_{ij}=\text{max}\{i+j-k,0\}$. 

$\A_{[2]}$ here is the part of $\V_{sl(2),q}$ built up from the
creation operators $x_{-1}$ only, {\it i.e.}\ the terms whose grade equals 
its degree. So its $(q,z)$-character is recovered
by substituting $uq^{-1}$ for $t$ in (4.1) and retaining only the constant
term  in $u$. We find that
$$ch_{\A_{[2]}(\V_{sl(2),k})}(q,z)=\sum_{m=0}^{2k}q^m\sum_{a=0}^{\text{min}\{
m,2k-m\}}(-1)^{m+a}\chi_{L(a)}(z)^2\ ,\eqno(4.2)$$
writing  $L(a)$ for the irreducible $a+1$-dimensional $sl(2)$-module, and hence 
$\A_{[2]}$ and $\A$ are isomorphic as $sl(2)$-modules.

\subhead IV.3. The root lattices\endsubhead
Consider first the $A_{N-1}$ root lattice, which can be identified with
the integer points $\vec{n}\in\bbZ^N$ with $\sum_i n_i=0$. Its automorphism
group is the symmetric group $Sym(N)$, together with $\vec{n}\mapsto-\vec{n}$,
so this will act on $\A_{[2]}$. Recall Section III.3. The `short' lattice
vectors $\vec{n}\in S_{A_{N-1}}$ are those whose coordinates $n_i$ all lie
in $\{\pm 1,0\}$; up to the $Sym(N)$ symmetry, we can take these to be
$\Lambda_\ell+\Lambda_{N-\ell}$, where $\ell\le \lfloor N/2\rfloor$ is the
number of
components equal to +1 and $\Lambda_i$ are the fundamental weights (the
natural basis for the dual lattice). There are $\left({N\atop 2\ell}\right)
\left({2\ell\atop \ell}\right)$ short vectors for a given $\ell$.

The number of  basis vectors (3.3) with $\alpha=0$ and grade $m$ is 
$\left({N\atop m}
\right)-\delta_{m,1}$, for a total (over all $m$) of $2^N-1$. This number for
$\alpha=\Lambda_\ell+\Lambda_{N-\ell}$ and grade $m$ is
$\left({N-2\ell\atop m}\right)$, for a total of $2^{N-2\ell}$. Therefore
the total dimension of the $C_2$-algebra is
$$\text{dim}\,\A_{[2]}(\V_{A_{N-1}})=2^N-1+\sum_{\ell=0}^{\lfloor N/2\rfloor}
2^{N-2\ell}\left({N\atop 2\ell}\right)\left({2\ell\atop \ell}\right)
=\left({2N\atop N}\right)-1\ .$$
By comparison, (2.6) tells us that Zhu's algebra has dimension
$$\text{dim}\,\A=\sum_{j=0}^{N-1}\left({N\atop j}\right)^2=\left({2N\atop N}
\right)-1\ .$$
So the $A_{N-1}$ root lattice is non-anomalous. (The `$-1$'s here, suggesting
a missing term, has an analogue in any affine $A$-series VOA, and is explained in
Section IV.5 below.)

The other root lattices can be handled similarly (in fact somewhat more
easily), with the result that only $E_8$ is anomalous. The short vectors for 
$E_8$ are 
0, a root, or the sum of 2 orthogonal roots. The $E_8$ Weyl group $W(E_8)$
acts transitively on each of those 3 sets, yielding 1-, 240-, and 2160-dimensional
$W(E_8)$-representations, respectively. $\A_{[2]}(\V_{E_8})$ is
the direct sum of the 2160-dimensional one, with 8 copies of the 240-dimensional
one, and 45 singlets, so is 4125-dimensional. But $\A_{[2]}$ also carries an
action of the $E_8$ Lie group (this is because the lattice VOA $\V_{E_8}$
is isomorphic to the affine algebra VOA $\V_{E_8,1}$), and in terms of this
it decomposes into $L(\Lambda_1)\oplus L(\Lambda_8)\oplus 2L(0)$, using
the node numbering conventions of Bourbaki/LiE (where $L(\Lambda_8)$ is
the 248-dimensional adjoint). By comparison, Zhu's algebra is 1-dimensional.

More generally, any (nontrivial) rational VOA with only 1 irreducible module
(these can be called {\it self-dual} VOAs) will be anomalous: Zhu's algebra
will be only 1-dimensional, because of (2.6), and $\A_{[2]}$ will always be
larger.

Incidentally, one-dimensional lattices  are easily shown to be 
non-ano\-ma\-lous.
Another simple fact: $\V_{L\oplus L'}=\V_L\otimes\V_{L'}$, so
$\V_{L\oplus L'}$ will be
anomalous iff either $\V_L$ or $\V_{L'}$ are.   Thus
it suffices to consider indecomposable lattices.

\subhead IV.4. Anomalous lattices\endsubhead
A lesson of the previous subsections is that among the most
accessible VOAs at least, the only anomalous
ones are anomalous for an elementary reason (namely, that they 
are self-dual). Because of this, it would be tempting to guess that 
anomalous VOAs are rare.

However, in this subsection and the next we shall give several VOAs which
are anomalous for subtle reasons. We suspect that in fact anomalous VOAs
are typical, for the following reason.

The paper [{\bf 2}] lists the indecomposable integral positive-definite 
lattices of
small dimension and determinant (the determinant will equal the number of 
irreducible $\V_L$-modules), 
and so can be regarded as providing some sort of random
sample of lattices. What we find is that, once we cross off from their
list root lattices and one-dimensional lattices,
which will automatically be non-anomalous, almost everything that remains
is anomalous!

Let $L$ be an $n$-dimensional even positive-definite lattice.
Using the analysis of Section III.3, we obtain the following (crude)
lower bound for the dimension of the $C_2$-algebra:
$$\text{dim}\,\A_{[2]}(\V_L)\ge \sum_{m=0}^{\mu+1}\left({n+m-1\atop m}
\right)+\left(n-\frac{1}{2}\right)M \ ,\eqno(4.3)$$
where $\mu$ is the minimum nonzero length-squared in $L$, and $M$ is the
number of lattice vectors with length-squared $\mu$. So the theta series
of $L$ starts like $1+Mq^{\mu/2}+\cdots$. To see (4.3), the sum over $m$
together with the term $-M/2$ bounds the number of vectors in (3.3) with
$\alpha=0$; each of the $M$ vectors with length-squared $\mu$ will also be
`small', and each of these will have at least $1+(n-1)$ vectors in (3.3).

There is no need to consider the determinant-1 lattices: they are all
anomalous. The 3 smallest indecomposable even lattices of determinant 2
are the root lattices $A_1$ and $E_7$, and the 15-dimensional lattice
called $D_{14}A_1[11]$. Consider the latter.
Its dim$\,\A$ is readily found to be $1^2+56^2=3137$.
It has $\mu=2$ and $M=366$, so (4.3) tells us its $\A_{[2]}$  is at least
6123-dimensional. Therefore $D_{14}A_1[11]$ is anomalous. This is typical
for the lattices collected in [{\bf 2}].

To get a clue as to what is special about the anomalous lattices, we should
ask what properties distinguish the $E_8$ root lattice from the other root
lattices. Of course, it is self-dual, but from our point of view this is
the wrong answer, as we now see there are plenty of nonself-dual anomalous lattices.
The most intriguing answer we have found is that $E_8$  is the only root lattice
whose {\it holes} do not lie in its dual. The holes of a lattice $L$ are
the points $\vec{x}$ in the ambient space $\bbR\otimes_\bbZ L$ whose distance
to any lattice point is a local maximum. If the hole is a global maximum,
it is called a {\it deep hole}. For example, $D_8=\{\vec{n}\in\bbZ^8\,|\,
\sum_i n_i\in 2\bbZ\}$ has deep holes at $(\frac{1}{2},\frac{1}{2},\ldots,
\pm\frac{1}{2})$ and a shallow hole at $(1,0,\ldots,0)$, and these all lie in
the weight lattice $D_8^*$. On the other hand, $E_8=\langle D_8,
(\frac{1}{2},\frac{1}{2},\ldots,\frac{1}{2})\rangle$ has a deep hole at
$(\frac{1}{2},\frac{1}{2},\ldots,-\frac{1}{2})$ and a shallow one at
$(\frac{1}{3},\frac{1}{3},\ldots,
\frac{1}{3},-\frac{1}{3})$, and of course neither lie in $E_8^*=E_8$.

Now, a vector in $\bbQ\otimes_\bbZ L$ (such as the holes of the integral 
lattice $L$)
will be a dual vector in some sublattice $L_0$ of $L$ of full dimension,
and thus will define an irreducible module of $\V_{L_0}$ and (lifting it
to $\V_L$) a (generally) twisted module of $\V_L$. The holes of a lattice $L$
should define special (perhaps twisted) $\V_L$-modules. For example,
the holes of $D_8$ all correspond to true $D_8$ modules, while the
holes of $E_8$ are twisted, coming from $D_8$ and $A_8$ sublattices (those twisted
modules have highest-weight spaces of dimension 16 and 9, respectively,
corresponding geometrically to the 16 and 9 vectors, respectively, of $E_8$ 
that are closest to the given hole). 
Recall the discussion at the end of Section II.3, where we explain that
$\A_{[2]}(\V)$ should see the twisted $\V$-modules $M$, in the sense that
there will be an embedding $M_0\otimes M_0^*\rightarrow \A_{[2]}^*$. The image
of this map may lie in the subspace of $\A_{[2]}$ coming
from the true $\V$-modules, but we would guess that
the twisted modules associated to holes would have an especially good chance
at landing in a new part of $\A_{[2]}$.

\subhead IV.5. Zhu's algebra and $C_2$-algebra for affine Lie algebras
\endsubhead
We understand the $A$-series quite well, at arbitrary rank and level,
with a conjectural description of $\A_{[2]}(\V_{sl(N),k})$ grade-by-grade
as an $sl(N)$-module. In particular, $\A_{[2]}$ at grade $m$
seems to be given by
$$
\A_{[2]}(\V)^{(m)}=\oplus_{\mu\in P_m^k}L(\mu)\otimes L(\mu)^*-\oplus_{\nu\in {P_m^k}'}
L(\nu)\otimes L(\nu)^*
\ , \eqno(4.4)
$$
where we define
$$\eqalignno{P_m^k&\,=\{\mu\in P_+\,|\,\mu_0\ge 0,\ t(\mu)\equiv m\
(\text{mod}\ N),\ t(\mu)\le m,\ N\mu_0+t(\mu)\ge m\}&(4.5)\cr
{P_m^k}'&=\{\nu\in P_+\,|\,\nu_0\ge 1,\ t(\nu)\equiv m-1\
(\text{mod}\ N),\ t(\nu)\le m-1,\ N\nu_0+t(\nu)\ge m\}\,,&}$$
 using $N$-ality $t(\mu)=\sum_{i=1}^{N-1} i\mu_i$, writing
$P_+$ for the $sl(N)$-weights with nonnegative Dynkin labels, and
setting {\it e.g.}\ $\mu_0=k-\sum_{i=1}^{N-1}\mu_i$.

This difference of modules appears because we are using $sl(N)$ rather
than $gl(N)$. The combinations $L(\mu)\otimes L(\mu)^*$ {\it etc.}\ arise 
ultimately because
of Peter-Weyl. Our conjecture is manifestly correct for grade $m\le k$, as the
null vector does not come in until $m=k+1$. If our
conjecture is correct, then the final nontrivial part of $\A_{[2]}$ will  
appear at grade $m=Nk$, where it will be a singlet.
It is easy to verify that our conjecture works for $sl(2)$, and
that our conjecture implies  $\V_{sl(N),k}$ {\it is not anomalous, for any
$N$ and $k$}. In fact, not only do the dimensions match, but $\A$ and $\A_{[2]}$
here are isomorphic as $sl(N)$-modules.

The few checks we have done {\it suggest} (although it is far too early to call
this even a conjecture) that likewise, $\V_{\g,k}$ is not anomalous
for any simple $\g$, except for $\g=E_8$. Of course $\V_{E_8,1}$ is isomorphic
to the self-dual lattice VOA $\V_{E_8}$, which being self-dual is anomalous
for elementary reasons. Remarkably, the $E_8$ VOAs are anomalous for all
levels except possibly $k=2$ [{\bf 11}].

Recall that both $\A(\V_{\g,k})$ and $\A_{[2]}(\V_{\g,k})$ carry an
 adjoint action of $\g$. For odd $k\ge 1$ the $E_8$-module $L(k\Lambda_1)$
 (again we follow the node numbering conventions of LiE/Bourbaki)
does not appear in Zhu's algebra as an irreducible summand, but appears in
the $C_2$-algebra. 
One can understand this in terms of $E_8$ twisted modules lifted
from $D_8$:
if we decompose the above module with respect to $D_8$ we get
$L(k\Lambda_1)^{e8} = L(2k\Lambda_1)^{d8} \oplus \cdots$.
Furthermore, none of the other $E_8$-modules that appear in Zhu's  
algebra can produce this $D_8$-module. On the other hand, in
$D_8$ we have
$$
L(k\Lambda_1)^{d8} \otimes L(k\Lambda_1)^{d8}
      = L(2k\Lambda_1)^{d8} \oplus\cdots \ .
$$
The module $L(k\Lambda_1)^{d8}$ is the highest-weight space of a  
level $k$ twisted $E_8$-module (restricted to $D_8$). This is
why $L(2k\Lambda_1)^{d8}$ must appear in the $C_2$-algebra of $E_8$, and 
implies
the $C_2$-algebra must be bigger than Zhu's   algebra (in fact it will
strictly contain it as an $E_8$-submodule).

For even levels $k>2$ the $E_8$-module $L((k-3)\Lambda_1+2\Lambda_2)$
is not in Zhu's algebra but appears in the $C_2$-algebra. However, we do  
not yet know how to obtain it from twisted modules.

We do not yet know whether level 2 is also anomalous. Curiously,
$E_8$ at level 2 has the only exceptional simple current ({\it i.e.}\ a simple
current not arising from an extended Dynkin diagram symmetry) among all the
affine algebras. 

\head V. Conjectures and questions\endhead

\roster

\item{{\it Clarify the role of holes in the lattice $L$ and $\A_{[2]}(\V_L)$.} 
We would guess that a lattice VOA $\V_L$ is anomalous whenever $L$ has a
hole not in its dual $L^*$.  Lattice VOAs are simple enough that we should
be able to completely characterise anomalous lattices.}

\item{{\it What is $\A_{[2]}(\V_{\g,k})$, grade by grade?}
In Section IV.5 we give a very satisfactory conjectural description of 
$\A_{[2]}(\V_{sl(N),k})$.
We have at present no idea what $\A_{[2]}(\V_{\g,k})$ looks like, grade by
grade, for the other simple $\g$.}

\item{{\it Clarify the relation between $\A_{[2]}$ and twisted modules.}
Do twisted
modules suffice to span $\A_{[2]}$? Can anything be said about how the images
of the $g$-twisted Zhu algebras $\A_g$ in $\A_{[2]}$ fit together, as
the automorphism $g$ varies?}

\item{{\it Cohomological interpretations of $\A_{[2]}^*/\A^*$.}}
See Section III.1.

\item{{\it The `gradation' of Zhu's algebra versus $C_2$-algebra.}
For the Lie algebra VOAs $\V_{\g,k}$, we give in Section III.4 an especially
clean description of the graded algebra associated to Zhu's algebra; this
should permit a direct comparison of it with $\A_{[2]}$ for these VOAs,
and perhaps a deeper understanding of $\A(\V_{\g,k})$ versus 
$\A_{[2]}(\V_{\g,k})$.}

\item{{\it Comparing related spaces.} 
Instead of considering the vacuum module $\V$, we can also
study the analogous question, {\it i.e.}\ whether 
$\text{dim}\,\A_{[2]}(M)=\text{dim}\,\A_{[1,1]}(M)$ for arbitrary modules $M$.
At least for the Virasoro minimal models with
$(p,q)=(5,2), (7,2), (9,2), (4,3), (5,3) ,(7,3)$ this seems to be the case
for all modules $M$. On the other hand,
the dimensions of {\it e.g.}\ $\A_{[3]}(\V)$ and $\A_{[1,1,1]}(\V)$ seem to 
already differ for the minimal models. (These calculations were performed
by Andy Neitzke.)  It seems that comparing $\A_{[2]}(\V)$ and
$\A(\V)$ is the most fundamental question here.}

\item{{\it Natural maps between $\A^*$ and $\A_{[2]}$?} The enveloping
algebra $U\g$ is a co-commutative Hopf algebra, and the polynomial algebra
$S\g$ is its Hopf dual. Of course the algebras $\A(\V_{\g,k})$ and
$\A_{[2]}(\V_{\g,k})$ are naturally quotients of $U\g$ and $S\g$, respectively.
Does something like this happen for general $\V$, and does this have any
significance?}

\endroster

\enddocument